\documentclass{elsart}
\usepackage{amsmath}
\usepackage{amssymb}
\journal{Mathematica Bohemica}
\date{}

\newcommand\tr{\mathrm{tr}}
\newcommand\A{{\mathcal A}}
\newcommand\U{{\mathcal U}}
\newcommand\B{{\mathcal B}}
\newcommand\la{l^2(\A)}
\def\tor{\mathrm{Tor}}

\def\gldim{\mathrm{gl.dim}}
\def\dirlim{\mathop{\varinjlim}\limits}

\begin{document}

\begin{frontmatter}

\title{Semisimplicity and Global Dimension of a Finite von Neumann Algebra}

\author{Lia Va\v s}

\address{Department of Mathematics, Physics and Computer Science,
University of the Sciences in Philadelphia, 600 S. 43rd St.,
Philadelphia, PA 19104}

\ead{l.vas@usip.edu}

\begin{abstract}
We prove that a finite von Neumann algebra ${\mathcal A}$ is
semi\-simple if the algebra of affiliated operators ${\mathcal U}$
of ${\mathcal A}$ is semisimple. When ${\mathcal A}$ is not
semisimple, we give the upper and lower bounds for the global
dimensions of ${\mathcal A}$ and ${\mathcal U}.$ This last result
requires the use of the Continuum Hypothesis.
\end{abstract}


\end{frontmatter}

\section{Introduction}

A finite von Neumann algebra $\A$ comes equipped with a faithful
and normal trace that enables us to define the dimension not just
of a finitely generated projective module over $\A$ but also of
any $\A$-module. This property makes $\A$ an interesting algebra
and creates the possibilities for various applications. For
applications of group von Neumann algebras in geometry and
algebra, see \cite{Lu5}.

The algebra $\A$ has some nice ring theoretic properties. It
mimics the ring $\Zset$ in such a way that every finitely
generated module is a direct sum of a torsion and torsion-free
part. The dimension faithfully measures the torsion-free part and
vanishes on the torsion part. Although not without zero-divisors
and, as we are going to see, rarely Noetherian, a finite von
Neumann algebra is $\Zset$-like: it is semihereditary (i.e., every
finitely generated submodule of a projective module is projective)
and has the classical quotient ring, constructed in the same way
as $\Qset$ is constructed from $\Zset.$ The classical ring of
quotients $\U$ of $\A$ can be defined solely within the operator
theory as the algebra of affiliated operators. Although $\U$ has
many nice properties as a ring, it is not necessarily semisimple
(Artinian and with trivial Jacobson radical) like $\Qset$ is.

Let us consider the conditions.
\begin{itemize}
\item[(1)] $\U$ is semisimple.
\item[(2)] $\A$ is $*$-isomorphic to the finite sum of algebras of $m_i\times m_i$
matrices over $L^{\infty}(n_i),$ $m_i>0,$ $n_i\geq 0,$ $i=1,\ldots, k$ for some $k>0.$
\item[(3)] $\A$ is isomorphic to the finite sum of rings of $m_i\times m_i$ matrices
over $\Cset^{n_i},$ $m_i>0,$ $n_i\geq 0,$ $i=1,\ldots, k$ for some $k>0.$
\item[(4)] $\A$ is semisimple.
\item[(5)] $\A$ has finite $\Cset$-dimension.
\end{itemize}

It is well known that the conditions (2) -- (5) are equivalent.
Also, it is not hard to see that conditions (2) -- (5) imply (1).
Here, we shall prove that (1) implies the rest of the conditions
(Theorem \ref{semisimple}).

If a ring is not semisimple, its global dimension measures how
close it is to being semisimple. The bounds for global dimension
of $\U$ and $\A$ will be given in the infinite dimensional case.
This result uses the Continuum Hypothesis (CH). Namely, in Theorem
\ref{GlavnaGlobalDim}, we shall show:
\begin{enumerate}
\item (CH) If $\dim_{\Cset}\A=\aleph_1$ then $\gldim \U=2$ and $2\leq \gldim \A\leq 3.$
\item (CH) If $\dim_{\Cset}\A=\aleph_n,$ $n>0,$ then $2\leq\gldim \U\leq n+1$ and $2\leq\gldim \A \leq n+2.$
\end{enumerate}

The paper is organized as follows. In Sections \ref{As} and
\ref{Us}, we list some results on a finite von Neumann algebra and
its algebra of affiliated operators. In Sections
\ref{AlgPreliminaries} and \ref{OTPreliminatires}, we list the preliminary
facts and results that we need. In Section \ref{Main}, we prove
the result on the semisimplicity. In Section \ref{GlobalDim}, we
give the upper and lower bounds for the global dimension of
non-semisimple $\A$ and $\U.$

The paper is written to be accessible both to an algebraist and an
operator theorist, so sometimes even well known results from the
fields are referenced for the sake of readability.

\section{Finite von Neumann Algebras} \label{As}

Let $H$ be a Hilbert space and $\B(H)$ be the algebra of bounded
operators on $H$. The space $\B(H)$ is equipped with five
different topologies: norm, strong, ultrastrong, weak and
ultraweak. The statements that a $*$-closed unital subalgebra $\A$
of ${\B}(H)$ is closed in weak, strong, ultraweak and ultrastrong
topologies are equivalent. For details see \cite{D} or
\cite{KadRing} (Theorem 5.3.1).

A {\em von Neumann algebra} $\A$ is a $*$-closed unital subalgebra
of ${\B}(H)$ which is closed with respect to weak (equivalently
strong, ultraweak, ultrastrong) operator topology.

A $*$-closed unital subalgebra $\A$ of ${\B}(H)$ is a von Neumann
algebra if and only if $\A= \A''$ where $\A'$ is the commutant of
$\A.$ The proof can be found in \cite{KadRing} (Theorem 5.3.1).

Let $Z(\A)$ denotes the center of $\A.$ A von Neumann algebra $\A$
is {\em finite} if there is a linear function $\tr_{\A}:
\A\rightarrow Z(\A)$ called {\em center-valued (or universal)
trace} uniquely determined by the properties that
\begin{enumerate}
\item $\tr_{\A}(ab) =\tr_{\A}(ba).$

\item $\tr_{\A}(a^*a)\geq 0.$

\item $\tr_{\A}$ is {\em normal}: it is continuous with respect to
ultraweak topology.

\item $\tr_{\A}$ is {\em faithful}: $\tr_{\A}(a)=0$ for some
$a\geq 0$ (i.e. $a = bb^*$ for some $b\in\A$) implies $a=0$.
\end{enumerate}

The trace function extends to matrices
over $\A$ in a natural way: the trace of a matrix is the sum of
the traces of the elements on the main diagonal. This provides us
with a way of defining a convenient notion of, not necessarily
integer valued, dimension of any module.

If $P$ is a finitely generated projective $\A$-module, there exist
$n$ and $f:\A ^n\rightarrow\A^n$ such that $f=f^2=f^*$ and the
image of $f$ is $P.$ Then, {\em the dimension} of $P$ is
\[ \dim_{\A}(P)=\tr_{\A}(f).\] For details see \cite{Lia2}. The center-valued dimension was also studied in \cite{Lu5}.

If $M$ is any $\A$-module, {\em the dimension} $\dim_{\A}(M)$ is
defined as \[\dim_{\A}(M)=\sup \{ \dim_{\A}(P) |\; P \mbox{
fin. gen. projective submodule of }M\}\]
where the supremum on the
right side is an element of $Z(\A)$ if it exists and is a new
symbol $\infty$ otherwise. We define
$a+\infty=\infty+a=\infty=\infty+\infty$ and $a\leq \infty$ for
every $a\in Z(\A).$

The dimension satisfies the following properties.
\begin{enumerate}
\item If $\;0\rightarrow M_0\rightarrow M_1\rightarrow M_2\rightarrow
0$ is a short exact sequence of $\A$-modules, then $
\dim_{\A}(M_1)= \dim_{\A}(M_0)+\dim_{\A}(M_2).$

\item If $M$ is a finitely generated projective module, then
$\;\dim_{\A}(M)=0$ if and only if $M=0.$
\end{enumerate}

The proof can be found in \cite{Lia2}. Part (1) follows from the first part of Proposition 13 in \cite{Lia2}. Part (2) follows from Theorem 17 in \cite{Lia2} (Additivity). For more properties of dimension, see \cite{Lu5} or \cite{Lia2}.

As a ring, a finite von Neumann algebra $\A$ is {\em
semihereditary} (i.e., every finitely generated submodule of a
projective module is projective or, equivalently, every finitely
generated ideal is projective). This follows from two facts.
First, every von Neumann algebra is an $AW^*$-algebra and, hence,
a Rickart $C^*$-algebra (see Chapter 1.4 in \cite{Be2}). Second, a
$C^*$-algebra is semihereditary as a ring if and only if it is
Rickart (see Corollary 3.7 in \cite{AG}). The fact that $\A$ is
Rickart also gives us that $\A$ is  {\em nonsingular} (see 7.6 (8)
and 7.48 in \cite{Lam}).

Note also that every statement about left ideals over $\A$ can be
converted to an analogous statement about right ideals. This is
the case because $\A$ is a ring with involution (which gives a
bijection between the lattices of left and right ideals and which
maps a left ideal generated by a projection to a right ideal
generated by the same projection).

\section{Algebras of Affiliated Operators} \label{Us}

A finite von Neumann algebra $\A$ is a pre-Hilbert space. Let $\la$
denote the Hilbert space completion of $\A.$ $\A$ can be
identified with the set of $\A$-equivariant bounded operators on
$\la,$ ${\B}(\la)^{\A},$ using the right regular representations
(see section 9.1.4 in \cite{Lu5} for details).

Let $a$ be a linear map $a: \mathrm{dom}\; a$ $\rightarrow$ $\la$,
$\mathrm{dom}\; a \subseteq\la.$ We say that $a$ is {\em
affiliated to $\A$} if
\begin{itemize}
\item[i)] $a$ is densely defined (the domain
$\mathrm{dom}\; a$ is a dense subset of $\la);$
\item[ii)] $a$ is closed (the graph of $a$ is closed in $\la \oplus \la);$
\item[iii)]  $ba = ab$ for every $b$ in the commutant of $\A.$
\end{itemize}
Let $\U=\U(\A)$ denote the {\em algebra of operators affiliated
to} $\A$.

\begin{prop}
Let $\A$ be a finite von Neumann algebra and $\U = \U(\A)$ its
algebra of affiliated operators.
\begin{enumerate}
\item $\A$ is an Ore ring and $\U$ is the classical ring of quotients
$Q_{\mathrm{cl}}(\A)$ of $\A.$
\item $\U$ is a von Neumann regular (fin. gen. submodule of fin.
gen.projective module is a direct summand), left and right
self-injective ring equal to the maximal ring of quotients
$Q_{\mathrm{max}}(\A).$
\item $\U$ is the injective envelope $E(\A)$ of $\A$.
\item The set of projections (idempotents) in $\U$ is the same as the set of projections (idempotents) in $\A.$
\end{enumerate}
\label{UisRingOfQuotients}
\end{prop}
The proof of 1. can be found in \cite{Lu5} (Theorem 8.22). The
proof of 2. is in \cite{Be1} (Lemma 1, Theorem 2, Theorem 3). 3.
follows from Theorem 13.36 in \cite{Lam} (note that $\A$ is
nonsingular).  4. follows from Theorem 1 in Section 48 and
Corollary 1 in Section 49 from \cite{Be2}. For a review of the ring theoretic notions of Ore ring, classical and maximal ring of quotients, self-injective ring and injective envelope see \cite{Lam}.

From this proposition it follows that the algebra $\U$ can be
defined using purely algebraic terms (ring of quotient, injective
envelope) on one hand and using just the language of operator
theory (affiliated operators) on the other.

$K_0(\A)$ and $K_0(\U)$ are isomorphic. The isomorphism
$\mu: K_0(\A)\cong K_0(\U)$ is induced by the map Proj($\A$)
$\rightarrow$ Proj($\U$)  given by $[P]\mapsto[\U\otimes_{\A} P]$
for any finitely generated projective module $P$ (Theorem 8.22 in
\cite{Lu5}). In \cite{Lia1} (Theorem 5.2) the explicit description of the map
Proj($\U$) $\rightarrow$ Proj($\A$) that induces the inverse of
the isomorphism $\mu$ is obtained. Namely, the following holds.

\begin{thm}
There is an one-to-one correspondence between direct summands of
$\A$ and direct summands of $\U$ given by $I \mapsto$
$\U\otimes_{\A}I = E(I).$ The inverse map is given by $L\mapsto
L\cap \A.$ This correspondence induces an isomorphism of monoids
$\mu :\mathrm{Proj}(\A)$ $\rightarrow$ $\mathrm{Proj}(\U)$ and an
isomorphism $\mu: K_0(\A)\rightarrow K_0(\U)$ given by
$[P]\mapsto[\U\otimes_{\A}P]$ with the inverse
$[Q]\mapsto[Q\cap\A^n]$ if $Q$ is a direct summand of $\U^n.$
\label{K0thm}
\end{thm}

For proof, see Theorem 5.2 in \cite{Lia1}.

\section{Algebraic Preliminaries}\label{AlgPreliminaries}

Let $R$ be a ring with unit. Let $M_n(R)$ denotes the ring of
$n\times n$ matrices over $R.$

\subsection{ }\label{WeddernburnArtin} Wedderburn-Artin Theorem
asserts that a ring $R$ is semisimple (Artinian with trivial
Jacobson radical) if and only if there are positive integers $m$
and $n_i,$ $i=1,\ldots,m$ and division rings $D_i,$ $i=1,\ldots,m$
such that $R$ is isomorphic to the product of matrix rings
$M_{n_i}(D_i),$ $i=1,\ldots,m.$ This result can be found in most
of the algebra textbooks (e.g. Theorem 3.3 in \cite{Hu}).

\subsection{ }\label{Morita} Morita Invariance Theorem asserts that $K_0(R)\cong
K_0(M_n(R))$ for every ring $R$ and every positive integer $n$
(see, for example, \cite{Rosenberg}, Theorem 1.2.4).

\subsection{ }\label{K0semisimple} If $R$ is a
semisimple ring, $K_0(R)$ is a finitely generated abelian group.
To show this, let $R$ be a semisimple ring. By Wedderburn-Artin
Theorem, $R$ is isomorphic to $\prod_{i=1}^m M_{n_i}(D_i)$ for
some $m>0, n_i>0,$ $i=1,2,\ldots,m$ and division rings $D_i$. Then
\[ K_0(R) \cong K_0(\prod_{i=1}^m M_{n_i}(D_i)) \cong
\bigoplus_{i=1}^m K_0( M_{n_i}(D_i))  \cong \bigoplus_{i=1}^m
K_0(D_i)  \cong \Zset^m.
\]
For details see 1.2.8  and 1.1.6 in \cite{Rosenberg}.

\subsection{ }\label{MinIdeals} A ring $R$ with unit is
semisimple if and only if R is a direct sum of minimal left ideals
of the form $Re_i,$ $i=1,\ldots,m$ where $e_1,e_2,\ldots e_m$ are
orthogonal idempotents with $e_1+e_2+\ldots+e_m=1.$ This result
can be found in various algebra textbooks (e.g. Theorem 3.7 in
\cite{Hu}).

Note that if $\A$ is a finite von Neumann algebra $\A,$ $p$ a
projection of $\A$ and the left ideal $\A p$ is minimal, then $p$
is a minimal projection. This is because for the projections
$p,q\in\A,$ $p=q$ if and only if $\A p=\A q$ (as in any ring with
involution). Since $\U$ does not have any new projections (part
(4) of Proposition \ref{UisRingOfQuotients}), the same holds for
$\U:$ if the left ideal $\U p$ is minimal, then $p$ is a minimal
projection.

Also note that if $e$ is an idempotent in $\A,$ then there is a
projection $p$ such that $\A e=\A p.$ This is true for every
Rickart $*$-ring (see, for example, section 3 in \cite{Be2}) so it
holds for $\A.$ Since the algebra of affiliated operators $\U$ is
also Rickart, the same holds for $\U$ and an idempotent $e\in \U$.
In that case, note that $p$ is in $\A$ (part (4) of Proposition
\ref{UisRingOfQuotients}).

\subsection{ }\label{MnSemisimple} If $R$ is a ring with unit,
then $R$ is semisimple if and only if the ring of $n\times n$
$R$-matrices $M_n(R)$ is semisimple. This follows from the
following two well known facts. First, $R$ is Artinian if and only
if $M_n(R)$ is Artinian (Exercise 5 in 8.1 of \cite{Hu}). Second,
$J(M_n(R))=M_n(J(R)),$ where $J(R)$ denotes the Jacobson radical
of the ring $R$ (Exercise 13 in 9.2 of \cite{Hu}).

Also, a product of rings $\prod_{i=1}^n R_i$ is semisimple if and only if each $R_i,$ $i=1,\ldots, n$ is semisimple (see Theorem 2.17 in 9.2 and Corollary 1.6 in 8.1 in \cite{Hu}).

\subsection{ }\label{MnEnvelope} If $E(R)$ is an injective
envelope of a ring $R$, it is easy to see that $M_n(E(R))$ is an
injective $M_n(R)$-module which is an essential extension of
$M_n(R).$ Thus, $M_n(E(R))=E(M_n(R)).$

This gives us that $\U(M_n(\A))=M_n(\U(\A))$ for a finite von
Neumann algebra $\A$ and its algebra of affiliated operators
$\U=E(\A).$

\section{Operator Theory Preliminaries}
\label{OTPreliminatires}

\subsection{ }\label{types} If $\A$ is a finite von Neumann
algebra, $\A$ has a natural and unique decomposition as a direct
sum of von Neumann algebras of types $I_n$ for positive
integers $n,$ and an algebra of type $II_1.$ For the definition of types of
von Neumann algebras and more details, see Section 6.5 in
\cite{KadRingtwo}. By definition of the type $I_1,$ $\A$ is
abelian if and only if $\A$ is of type $I_1.$

\subsection{ }\label{typeIn}
If $\A$ is of type $I_n$, $\A$ is $*$-isomorphic to the algebra
$M_n(Z(\A))$ where $Z(\A)$ is the center of $\A$ (Theorem 6.6.5.
in \cite{KadRingtwo}).

\subsection{ }\label{typeII1}  If a finite von Neumann algebra
is of type $II_1,$ the group $K_0(\A)$ is isomorphic to the group
$Z(\A)^{\Zset_2} = \{\;a\in Z(A)\;|\;a = a^*\;\},$ the subgroup of
$\Zset_2$-invariants of $Z(\A)$ with the action of $\Zset_2$ by
involution. This follows from the proof of Theorem 8.4.4. in
\cite{KadRingtwo}. Also, see Theorem 9.13. in \cite{Lu5}.

As a consequence of this, $K_0(\A)$ of a nontrivial von Neumann
algebra of type $II_1$ is not finitely generated. Indeed, if $\A$
is nontrivial, the group $Z(\A)^{\Zset_2}$ is nontrivial as well.
If nontrivial, the group $Z(\A)^{\Zset_2}$ contains a copy of
$\Rset$ because it contains all the projections of the form $r 1$
where $r\in \Rset$ and 1 is the unit of $\A.$

\subsection{ }
\label{LInfExamples}

If $S$ is a compact space and
$\mu$ a finite measure on the Borel algebra of $S,$ we use the standard notation
$L_2(S, \mu)$ and $L^{\infty}(S, \mu)$ in their usual sense. If $S$ is $\{1, 2,..., n\},$ $\{1,2,...\}=\aleph_0$
or $[0,1],$ the corresponding algebras are denoted by $L^{\infty}(n),$ $L^{\infty}(\aleph_0)$
and $L^{\infty}([0,1])$ respectfully. Note that the first one is semisimple and finitely generated while the other two are not.

\subsection{ }\label{abelian}
We use the term maximal abelian von Neumann algebra in its usual sense (e.g. see
\cite{KadRing}). Recall that
every abelian von Neumann algebra is $*$-isomorphic to a maximal
abe\-lian algebra. For the proof of this see section 9.4. in
\cite{KadRingtwo} or section I 7 in \cite{D}.

Every two $*$-isomorphic maximal abelian von Neumann algebras are
unitarily equivalent. Moreover, the Hilbert spaces on which these
two algebras act are isomorphic (Theorem 9.3.1.
\cite{KadRingtwo}).

\subsection{ }\label{MaxAbelian}
Theorem 9.4.1. from \cite{KadRingtwo}, asserts that every maximal
abelian von Neumann algebra that acts on a separable Hilbert space
is unitarily equivalent to exactly one of the algebras
$L^{\infty}(n),$ $L^{\infty}(\aleph_0),$ $L^{\infty} ([0,1]),$
$L^{\infty}(n)\oplus L^{\infty}([0,1])$ or
$L^{\infty}(\aleph_0)\oplus L^{\infty}([0,1]).$

\section{Semisimplicity}
\label{Main}

For the main result, we need the following lemma.

\begin{lem} If $\A$ is an abelian von Neumann algebra such that
the algebra of affiliated operators $\U$ is semisimple, then $\A$
is finite dimensional. \label{AbelianSemisimple}
\end{lem}
\begin{pf}
Let $\A$ be an abelian von Neumann algebra with semisimple $\U.$
There are orthogonal idempotents $e_i,$ $i=1,2,\ldots,n,$ such
that the left ideals $\U e_i$ are minimal $i=1,2,\ldots,n,$
$\sum_{i=1}^n e_i =1$ and $\U=\bigoplus_{i=1}^n \U e_i$ by
\ref{MinIdeals}. There are projections $p_i\in\A,$ such that $\U
e_i=\U p_i,$ $i=1,2,\ldots,n$ (\ref{MinIdeals}). The projections
$p_i,$ $i=1,2,\ldots,n,$ are minimal projections since the ideals
$\U p_i,$ $i=1,2,\ldots,n,$ are minimal (\ref{MinIdeals}).

$\A p_i$ is a direct summand of $\A$ and so $\A p_i=\A \cap
(\U\otimes_{\A}\A p_i)$ (by Theorem \ref{K0thm}). But
$\U\otimes_{\A}\A p_i=\U p_i$ and so $\A p_i= A\cap \U p_i.$

If $j\neq i,$ then $\A p_i\cap \A p_j =\A\cap(\U p_i\cap\U p_j)=
\A\cap(\U e_i\cap\U e_j)=0$ and so $\bigoplus_{i=1}^n \A
p_i\subseteq \A.$ Moreover, \[\A=\A (\sum_{i=1}^n e_i)\subseteq
\sum_{i=1}^n \A e_i\subseteq \sum_{i=1}^n \A\cap \U e_i =
\sum_{i=1}^n \A \cap \U p_i = \sum_{i=1}^n \A p_i =
\bigoplus_{i=1}^n \A p_i.\]

The spaces  $\A p_i,$ $i=1,2,\ldots,n,$ are one-dimensional so $\A$ is finite dimensional.
\end{pf}

\begin{thm} Let $\A$ be a finite von Neumann algebra with the
algebra of affiliated operators $\U$. The following are
equivalent:
\begin{enumerate}
\item $\U$ is semisimple.
\item $\A$ is $*$-isomorphic to the finite sum of algebras of $m_i\times m_i$
matrices over $L^{\infty}(n_i),$ $m_i>0,$ $n_i\geq 0,$ $i=1,\ldots, k$ for some $k>0.$
\item $\A$ is isomorphic to the finite sum of rings of $m_i\times m_i$ matrices
over $\Cset^{n_i},$ $m_i>0,$ $n_i\geq 0,$ $i=1,\ldots, k$ for some $k>0.$
\item $\A$ is semisimple.
\item $\A$ has finite $\Cset$-dimension.
\end{enumerate}
\label{semisimple}
\end{thm}
It is well known that the conditions (2) -- (5) are equivalent.
Also, it is not hard to see that conditions (2) -- (5) imply (1).
The main result here is that (1) implies the rest of the
conditions.
\begin{pf}
$(1)\Rightarrow(2).$ Let $\U$ be semisimple. Then $K_0(\U)$ is
finitely generated by \ref{K0semisimple}. Since $K_0(\U)\cong
K_0(\A)$ (Theorem \ref{K0thm}), $K_0(\A)$ is finitely generated as well.
Hence, $\A$ has no summand of type $II_1$ and has just finitely many summands $I_{m_i},$
$m_i>0,$ $i=1,\ldots, k$ by \ref{types}
and \ref{typeII1}. Let us denote these summands by $\A_{m_i}.$
For each $i,$ $\A_{m_i}$ is $*$-isomorphic to the
algebra $M_{m_i}(Z(\A_{m_i}))$
by \ref{typeIn}. Note that $K_0(Z(\A_{m_i}))=K_0(\A_{m_i})$ by
\ref{Morita} and $\U(\A_{m_i}) = \U(M_m(Z(\A_{m_i}))) = M_m(\U(Z(\A_{m_i})))$ by
\ref{MnEnvelope}. By \ref{MnSemisimple},
$\U(Z(\A_{m_i}))$ is
semisimple. Thus, $Z(\A_{m_i})$ is finite dimensional by Lemma \ref{AbelianSemisimple}.

$Z(\A_{m_i})$ is $*$-isomorphic to a maximal abelian algebra by
\ref{abelian}. As $Z(\A_{m_i})$ is finite dimensional, it is $*$-isomorphic to $L^{\infty}(n_i)$ for
some nonnegative integer $n_i$ by \ref{MaxAbelian}.
Thus, $\A_{m_i}$ is $*$-isomorphic to the algebra of $m_i\times m_i$
matrices over $L^{\infty}(n_i).$

$(2)\Rightarrow(3).$ $L^{\infty}(n)$ is isomorphic to $\Cset^n$ as
rings.

$(3)\Rightarrow(4).$ If condition (3) is satisfied, then $\A$ is semisimple by \ref{WeddernburnArtin}
and \ref{MnSemisimple}.

$(4)\Rightarrow(1).$ If $\A$ is semisimple, then $\A$ is
self-injective and, hence, equal to its injective envelope $\U.$
Thus, $\U=\A$ is also semisimple.

$(5)\Rightarrow(4).$ If $\A$ has finite $\Cset$-dimension, then
$\A$ is $*$-isomorphic to the direct sum of finitely many algebras
of the form $M_{n_i}(\Cset),$ $i=1,\ldots,k.$ This follows from
Proposition 6.6.6, Theorem 6.6.1 and comments preceding
Proposition 6.6.6 from \cite{KadRingtwo} Thus, $\A$ is semisimple
by \ref{WeddernburnArtin}.

$(2)\Rightarrow(5)$ is clear.
\end{pf}
Clearly, $\A$ is abelian if and only if $m=1$ in (2) and (3).

\begin{cor} Let $\A$ be a finite von Neumann algebra with the
algebra of affiliated operators $\U$. The following conditions are
equivalent to conditions (1) -- (5) from Theorem \ref{semisimple}:
\begin{itemize}
\item[(6)] $\A$ is Noetherian (every ideal is finitely generated).
\item[(7)] $\U$ is Noetherian.
\item[(8)] $\A$ is hereditary (every ideal is projective).
\item[(9)] $\U$ is hereditary.
\item[(10)] $\A$ has finite universal dimension (every direct sum
of nonzero submodules contained in $\A$ is finite).
\item[(11)] $\U$ has finite universal dimension.
\end{itemize}
\label{CorollarySemisimple}
\end{cor}
Note that the only nontrivial part is that (8) implies any of the
other conditions.

\begin{pf} Clearly, (1) implies (7), (9) and (11), and (4) implies
(6), (8), and (10). Since (1) and (4) are equivalent, to show the
equivalence of all of the conditions, it is sufficient to prove
that each of the conditions (6) -- (11) implies (1). We shall show
$(6)\Rightarrow (7)\Rightarrow (1),$ $(8)\Rightarrow
(9)\Rightarrow(1)$ and $(11)\Leftrightarrow (10)\Rightarrow (1).$

$(6)\Rightarrow (7)\Rightarrow (1).$ The classical ring of
quotients of an Ore and Noetherian ring is Noetherian (Proposition
10.32 (6) in \cite{Lam}). Since $\U=Q_{\mathrm{cl}}(\A)$ and $\A$
is Ore, $(6)\Rightarrow(7).$ A von Neumann regular and Noetherian
ring is semisimple (Corollary 5.60 and Example 5.62a in
\cite{Lam}). Thus, $(7)\Rightarrow(1).$

$(8)\Rightarrow (9)\Rightarrow(1).$ A self-injective and
hereditary ring is semisimple (Theorem 7.52 in \cite{Lam}). Since
$\U$ is self-injective, $(9)\Rightarrow(1).$ We shall show
$(8)\Rightarrow (9),$ in two steps. First, we show that for every
ideal $J$ of $\U,$ $J=\U\otimes_{\A}(J\cap \A).$ Second, we show
that if $\A$ is hereditary, $\U\otimes_{\A}(J\cap \A)$ is
projective. This will give us that every ideal $J$ of $\U$ is
projective.

Let $J$ be an ideal of $\U.$ Consider first the case when $J$ is
finitely generated. Since $\U$ is semihereditary, $J$ is (finitely
generated and) projective. Thus, there is a positive integer $n$
such that $J$ is a direct summand of $\U^n.$ $\U$ is
self-injective, so $J$ is a direct summand of an injective module
and therefore $J$ is injective. So, the inclusion
$J\hookrightarrow\U$ splits, so $J$ is a direct summand of $\U.$
Then, Theorem \ref{K0thm} gives us that $J=\U\otimes_{\A}(J\cap
\A).$

If $J$ is any ideal of $\U,$ $J$ is the directed union of its
finitely generated submodules $J_i$ (directed with respect to the
inclusion maps). Then,
\[
\begin{array}{rcl}
\U\otimes_{\A}(J\cap \A)) & = & \U\otimes_{\A}((\dirlim J_i)\cap
\A))=\U\otimes_{\A}\dirlim (J_i\cap \A))= \\ & = & \dirlim
\U\otimes_{\A}(J_i\cap \A)) =\dirlim J_i=J.
\end{array}
\]

Now, let us show that $\U\otimes_{\A}(J\cap \A))=J$ is projective
for $\A$ hereditary. If $\A$ is hereditary, the module $J\cap \A$
is projective and so a direct summand of some free module
$\bigoplus\A$. But then, $\U\otimes_{\A}(J\cap \A))$ is a direct
summand of $\bigoplus \U$ and so, projective. This gives us that
every ideal of $\U$ is projective, so $\U$ is hereditary.

$(11)\Leftrightarrow (10)\Rightarrow (1).$ If a ring $R$ is Ore,
the uniform dimension of $R$ is equal to the uniform dimension of
$Q_{\mathrm{cl}}(R)$ (Corollary 10.35 in \cite{Lam}). Thus,
$(10)\Leftrightarrow (11).$ If a nonsingular ring $R$ has finite
universal dimension, then $Q_{\mathrm{max}}(R)$ is semisimple
(Theorem 13.40 in \cite{Lam}). Since $\A$ is nonsingular and
$\U=Q_{\mathrm{max}}(\A),$ $(10)\Rightarrow(1)$ follows.
\end{pf}

\section{Global dimensions of $\U$ and $\A$}\label{GlobalDim}

In this section, we shall examine the global dimension of rings
$\U$ and $\A$. The global dimension of a ring measures how close
the modules over that ring are to being projective, therefore how
close the ring is to being semisimple. The bounds for global
dimension of $\U$ and $\A$ will be given.

The global dimension of a ring $R$ is defined via the {\em
projective dimension} of a left $R$-module $M$ \[\mbox{pd}_R(M) =
\min\{\;n\;|\;M\mbox{ has a projective resolution of length
}n\;\}.\] If this minimum does not exist, we define pd$_R(M)$ to
be $\infty.$ Clearly, a left module $M$ is projective if and only
if pd$(M)$=0.

The {\em left global dimension} of a ring $R$ is
\[\mbox{l.gl.}\dim R = \sup\{\;\mbox{pd}_R(M)\;|\;M\mbox{ is a left }R\mbox{-module}\;\}.\]

The left global dimension can be computed using ideals solely:
\[\mathrm{l.gl.dim} R = \sup\{\;\mathrm{pd}(R/I)\;|\;I\mbox{ is a left }R\mbox{-ideal}\;\}.\]
See Corollary 5.51 in \cite{Lam} for details.

The right global dimension is defined similarly. If left and right
global dimensions of a ring are equal, we write just $\gldim R$
for l.gl.$\dim R = $ r.gl.$\dim R.$ This is the case for $\A$ and
$\U$ since they are rings with involution so every statement about
left ideals can be converted to an analogous statement about right
ideals.

Clearly, a ring $R$ is semisimple iff r.gl.$\dim R = 0$ iff
l.gl.$\dim R = 0$. Also, $R$ is left hereditary (every submodule
of a projective left module is projective) if and only if
l.gl.$\dim R\leq 1.$

We have seen that $\gldim \U = 0$ just if $\A$ is finite
dimensional. Suppose that  $\gldim \U = 1$. Then, $\U$ is
hereditary. But every self-injective and hereditary ring is
semisimple (we already used this in $(9)\Rightarrow(1)$ of
Corollary \ref{CorollarySemisimple}), so $\gldim \U = 0.$ So, if
$\A$ is infinite dimensional, the global dimension of $\U$ is at
least 2.

Tor functor defines another dimension of a ring. For review of Let $R$ be a ring
and $M$ a left $R$-module. The {\em weak dimension} of $M$ is
\[\mbox{wd}(M) = \sup\{\;n\;|\;\tor_n^R(\underline{\hskip0.3cm},
M)\neq 0\;\}.\] Clearly, $M$ is a flat left module if and only if
wd $(M)=0.$ If $M$ is a right module, we can define its weak
dimension as the supremum of dimensions $n$ of nonvanishing
$\tor^R_n(M, \underline{\hskip0.3cm}).$ It can be shown that the
supremum of weak dimensions of left modules is the same as the
supremum of weak dimensions of right modules and that is the same
as $\sup\{\;n\;|\;\tor^R_n(\underline{\hskip0.3cm},
\underline{\hskip0.3cm})\neq 0\;\}$ so, we can define the {\em
weak global dimension} of $R$ as \[\mbox{wd} R =
\sup\{\;n\;|\;\tor^R_n(\underline{\hskip0.3cm},
\underline{\hskip0.3cm})\neq 0\;\}.\] Since this definition is
left-right symmetric, we do not have to distinguish left and right
weak global dimension. For more details, see section 5D in \cite{Lam}.

A ring $R$ is von Neumann regular if and only if all modules are
flat (Theorem 4.21, \cite{Lam}). Thus, $R$ is von Neumann regular
if and only if wd $R = 0$. So,
\[\mathrm{wd} \U =0.\]

For any ring $R$, wd $R \leq 1$ if and only if a submodule of a
flat module is flat. Since all semihereditary rings have this
property (see Theorem 4.67 in \cite{Lam}), \[\mathrm{wd} \A\leq
1.\] There are von Neumann algebras with weak global dimension 1
(Example 2.9 in \cite{Lu5}).

The following theorem of Jensen (Theorem 5.2 in \cite{Jen2})
connects the global dimension of a ring with its cardinality and
its weak global dimension. Recall that $\aleph_0$ denotes the
first infinite cardinal, the cardinality of the set of integers.
Then, $\aleph_{n+1}$ is defined as the successor cardinal of
$\aleph_n,$ the least cardinal strictly larger than $\aleph_n.$

\begin{thm} If $R$ is a ring of cardinality $\aleph_n,$ $n\geq 0,$ then
\[\begin{array}{rcl}
\mathrm{l.gl.dim} R & \leq & \mathrm{wd
} R + n +1\mbox{ and}\\
\mathrm{ r.gl.dim }R & \leq & \mathrm{ wd }R + n +1.
\end{array}\] \label{UpperBounds}
\end{thm}

If $\A$ is a finite von Neumann algebra, the cardinality of $\A$
is at least the continuum $c$ (the cardinality of $\Cset$) since
$\A$ contains a copy of the set of complex numbers $\{\;z
1_{\A}\;|\;z\in \Cset\;\}$ where $1_{\A}$ is the identity operator
in $\A$. Also, since $\U=Q_{\mathrm{cl}}(\A),$
$\A\subseteq\U\subseteq\A\times\A.$ The cardinality of $\A$ is the
same as the cardinality of $\A\times\A$ since both are infinite,
so the cardinality of $\A$ and $\U$ are the same.

Since wd $\U=0$ and wd $\A\leq 1,$ the Theorem \ref{UpperBounds}
gives us the following.
\begin{cor} If $\A$ is a finite von Neumann algebra of cardinality $\aleph_n,$ $n>0,$ with the
algebra of affiliated operators $\U,$ then
\[\begin{array}{rcl}
\gldim \U & \leq & n +1\mbox{ and}\\
\gldim \A & \leq & n +2.
\end{array}\] \label{bounds}
\end{cor}

If we were to use this result, we would like to identify the
cardinality of $\A$ as one of $\aleph$'s. Note that the case of
$\aleph_{\lambda}$ when lambda is not finite in Jensen's theorem
is trivial.

If $V$ is an infinite dimensional complex space, its cardinality
and its dimension over $\Cset$ are closely connected. Namely, if
$\dim_{\Cset} V=\lambda,$ where $\lambda$ is infinite, then the
cardinality of $V$ is equal to the cardinality of finite subsets
of the set $\Cset \times \lambda.$ This cardinality is the same as
cardinality of $\Cset \times \lambda.$ Since $\lambda$ is
infinite, this is the maximum of the continuum $c$ and $\lambda.$

Thus, using Corollary \ref{bounds} for infinite dimensional $\A$
implies identifying the maximum of the dimension of $\A$ and the
continuum $c$ as one of the $\aleph$'s. That requires the use of
the {\em Continuum Hypothesis} (CH). Recall that CH states that
\[ c = \aleph_1.\]

In the sequel, we shall emphasize the use of CH.

\begin{thm} Let $\A$ be a finite von Neumann algebra with the
algebra of affiliated operators $\U,$ then
\begin{enumerate}
\item $\dim_{\Cset}\A<\aleph_0$ if and only if $\gldim \U=\gldim \A=0.$
\item (CH) If $\dim_{\Cset}\A=\aleph_1$ then $\gldim \U=2$ and $2\leq \gldim \A\leq 3.$
\item (CH) If $\dim_{\Cset}\A=\aleph_n,$ $n>0,$ then $2\leq\gldim \U\leq n+1$ and $2\leq\gldim \A \leq n+2.$
\end{enumerate}
\label{GlavnaGlobalDim}
\end{thm}
\begin{pf}
(1) is proven in Theorem \ref{semisimple}.

If $\dim_{\Cset}\A=\aleph_1$ (note that $\dim_{\Cset}\A$ cannot be $\aleph_0$ as $\A$ is a Banach space) then the cardinalities
of $\A$ and $\U$ are both $\aleph_1.$ So, $\gldim \U\leq 2$ by
Corollary \ref{bounds}. If $\gldim \U\leq 1,$ then $\U$ is
hereditary and therefore semisimple (Corollary
\ref{CorollarySemisimple}). But then $\dim_{\Cset}\A<\aleph_0.$
Thus, $\gldim \U=2.$

Similarly, $\gldim \A\leq 3$ by Corollary \ref{bounds}. But
$\gldim \A\leq 1$ ($\A$ hereditary) is equivalent with the
conditions from Theorem \ref{semisimple} and Corollary
\ref{CorollarySemisimple} and implies $\dim_{\Cset}\A<\aleph_0.$
So, $\gldim \A\geq 2.$

(3) is proven analogously.
\end{pf}

Let D denote the statement: if $\dim_{\Cset}\A=\aleph_1,$ then
$\gldim \U=2.$ We have seen that CH implies D. The following
questions are open:
\begin{enumerate}
\item Does D hold without assuming CH? If the answer is yes, the proof would probably be very enlightening.
If the answer is no, the next question will be of great interest.
\item Does D imply CH? In other words, is D
equivalent with CH? If so, we have another equivalent of CH in our
hands.
\item Can the bounds for global dimension of $\A$ and $\U$ be narrowed i.e. can Theorem \ref{GlavnaGlobalDim} be
improved?
\item What can we say about global dimensions of $\A$ and $\U$ if
the $\Cset$-dimension of $\A$ is $\aleph_{\lambda}$ with $\lambda
$ an infinite ordinal?
\end{enumerate}

\section*{Acknowledgments}

The author is grateful to Peter Linnell for pointing out some incorrect statements in an earlier version of the text. The author would also like to thank the referee for insightful suggestions.

\bibliographystyle{amsalpha}

\end{document}